# A note on the Lindeberg condition for convergence to stable laws in Mallows distance

EURO G. BARBOSA[1] and CHANG C.Y. DOREA[2]

[1]*Banco Central do Brasil, SBS Quadra 3, 70074-900 Brasília-DF, Brazil.*
*E-mail: euro.barbosa@bcb.gov.br*
[2]*Departamento de Matemática, Universidade de Brasilia, 70910-900 Brasilia-DF, Brazil.*
*E-mail: changdorea@unb.br*

We correct a condition in a result of Johnson and Samworth (*Bernoulli* **11** (2005) 829–845) concerning convergence to stable laws in Mallows distance. We also give an improved version of this result, setting it in the more familiar context of a Lindeberg-like condition.

Theorem 5.2 of [1] considers a fixed parameter $\alpha \in (0,2)$, an independent sequence of random variables $X_1, X_2, \ldots$ with $S_n = (X_1 + \cdots + X_n)/n^{1/\alpha}$ and a random variable $Y$ with an $\alpha$-stable distribution. Theorem 5.2 claims that if there exist (independent) copies $Y_1, Y_2, \ldots$ of $Y$ satisfying

$$\frac{1}{n} \sum_{i=1}^{n} \mathbb{E}\{|X_i - Y_i|^\alpha \mathbb{1}(|X_i - Y_i| > b)\} \to 0 \tag{1}$$

as $b \to \infty$, then $S_n$ (possibly shifted) converges to $Y$ in Mallows distance $d_\alpha$. The proof given for Theorem 5.2 requires simultaneous control of $b$ and $n$, which is not provided by (1) as stated. Although the result could be corrected by adding "$\sup_n$" to the beginning of (1) and with other small modifications, we instead provide a more natural Lindeberg condition. We also change the centering, providing explicit expressions for the centering sequence for the case $\alpha \in (1,2)$. This is, in fact, a coupling theorem. Indeed, for $\alpha \in [1,2)$, if the Mallows distance between the distributions $F_X$ and $F_Y$ of $X$ and $Y$ is finite, then the random variables $X$ and $Y$ are highly dependent, in the sense that $d_\alpha^\alpha(X, Y) = \mathbb{E}|X - Y|^\alpha$ provided the joint distribution of $(X, Y)$ is $F_X \wedge F_Y$.

**Theorem 1.** *Fix $0 < \alpha < 2$. Let $(X_1, Y_1), (X_2, Y_2), \ldots$ be a sequence of independent pairs such that $Y_1, Y_2, \ldots$ are copies of an $\alpha$-stable random variable $Y$, and such that for all*

$b > 0$, *we have*

$$\lim_{n\to\infty} \frac{1}{n} \sum_{i=1}^{n} \mathbb{E}\{|X_i - Y_i|^\alpha \mathbb{1}(|X_i - Y_i| > bn^{(2-\alpha)/2\alpha})\} = 0. \tag{2}$$

*Then, writing $S_n = (X_1 + \cdots + X_n)/n^{1/\alpha}$, there exists a sequence of constants $(c_n)$ such that $\lim_{n\to\infty} d_\alpha(S_n - c_n, Y) = 0$. Moreover, when $\alpha \in (1, 2)$, we may take $c_n = n^{-1/\alpha} \times \sum_{i=1}^{n} \mathbb{E}X_i - \mathbb{E}Y$.*

**Proof.** By Corollary 1.2.9 of [2],

$$\frac{1}{n^{1/\alpha}} \sum_{i=1}^{n} Y_i \stackrel{d}{=} \begin{cases} Y + \mu n^{1-1/\alpha} - \mu, & \text{if } \alpha \neq 1, \\ Y + \frac{2}{\pi}\sigma\beta \log n, & \text{if } \alpha = 1. \end{cases}$$

Here, the constants $\mu \in \mathbb{R}$, $\sigma \geq 0$ and $\beta \in [-1, 1]$ are, respectively, the shift, scale and skewness parameters of the stable law of $Y$ (see, e.g., [2], page 5), so for $\alpha \in (1, 2)$, we may take $\mu = \mathbb{E}Y$. We first treat the case $\alpha \in (1, 2)$. With $c_n$ defined as in the statement of the theorem,

$$S_n - c_n - Y \stackrel{d}{=} n^{-1/\alpha} \sum_{i=1}^{n} (U_i - \mathbb{E}U_i + V_i - \mathbb{E}V_i),$$

where, writing $\delta = \frac{2-\alpha}{2\alpha}$,

$$U_i = (X_i - Y_i)\mathbb{1}(|X_i - Y_i| \leq bn^\delta),$$
$$V_i = (X_i - Y_i)\mathbb{1}(|X_i - Y_i| > bn^\delta).$$

Using Lyapunov's inequality and the fact that $|U_i| \leq bn^\delta$, we have

$$\mathbb{E}\left\{\left|\sum_{i=1}^{n}(U_i - \mathbb{E}U_i)\right|^\alpha\right\} \leq \left[\mathbb{E}\left\{\left|\sum_{i=1}^{n}(U_i - \mathbb{E}U_i)\right|^2\right\}\right]^{\alpha/2} = \left(\sum_{i=1}^{n} \operatorname{Var} U_i\right)^{\alpha/2}$$
$$\leq b^\alpha n^{(1+2\delta)\alpha/2} = b^\alpha n. \tag{3}$$

Similarly, a von Bahr–Esseen moment bound given as equation (12) in [1] yields

$$\mathbb{E}\left\{\left|\sum_{i=1}^{n}(V_i - \mathbb{E}V_i)\right|^\alpha\right\} \leq 2\sum_{i=1}^{n}\mathbb{E}(|V_i - \mathbb{E}V_i|^\alpha) \leq 2^{\alpha+1}\sum_{i=1}^{n}\mathbb{E}(|V_i|^\alpha). \tag{4}$$

Thus, by (3) and (4), we find that for $\alpha \in (1, 2)$,

$$d_\alpha^\alpha(S_n - c_n, Y) \leq \mathbb{E}\{|S_n - c_n - Y|^\alpha\}$$

$$\leq \frac{2^{\alpha-1}}{n}\mathbb{E}\left\{\left|\sum_{i=1}^{n}(U_i - \mathbb{E}U_i)\right|^{\alpha}\right\} + \frac{2^{\alpha-1}}{n}\mathbb{E}\left\{\left|\sum_{i=1}^{n}(V_i - \mathbb{E}V_i)\right|^{\alpha}\right\}$$

$$\leq 2^{\alpha-1}b^{\alpha} + \frac{2^{2\alpha}}{n}\sum_{i=1}^{n}\mathbb{E}\{|X_i - Y_i|^{\alpha}\mathbb{1}(|X_i - Y_i| > bn^{\delta})\}.$$

We deduce from condition (2) that $\limsup_{n\to\infty} d_{\alpha}^{\alpha}(S_n - c_n, Y) \leq 2^{\alpha-1}b^{\alpha}$. However, $b > 0$ was arbitrary, so the result follows.

When $\alpha \in (0, 1]$ and condition (2) holds, we can find a sequence $(b_n)$ converging to zero with

$$\lim_{n\to\infty}\frac{1}{n}\sum_{i=1}^{n}\mathbb{E}\{|X_i - Y_i|^{\alpha}\mathbb{1}(|X_i - Y_i| > b_n n^{(2-\alpha)/2\alpha})\} = 0.$$

In this case, we should define

$$c_n = \begin{cases} n^{-1/\alpha}\sum_{i=1}^{n}\mathbb{E}\{(X_i - Y_i)\mathbb{1}(|X_i - Y_i| \leq b_n n^{\delta})\} + \mu n^{1-1/\alpha} - \mu, & \text{for } 0 < \alpha < 1, \\ n^{-1/\alpha}\sum_{i=1}^{n}\mathbb{E}\{(X_i - Y_i)\mathbb{1}(|X_i - Y_i| \leq b_n n^{\delta})\} + \frac{2}{\pi}\sigma\beta \log n, & \text{for } \alpha = 1. \end{cases}$$

Then, with the same definitions of $U_i$ and $V_i$, except with $b$ replaced by $b_n$, we have

$$S_n - c_n - Y \stackrel{d}{=} n^{-1/\alpha}\sum_{i=1}^{n}(U_i - \mathbb{E}U_i + V_i).$$

The argument now mimics the case $\alpha \in (1, 2)$. Using analogues of the bounds (3) and (4), we find

$$d_{\alpha}^{\alpha}(S_n - c_n, Y) \leq b_n^{\alpha} + \frac{1}{n}\sum_{i=1}^{n}\mathbb{E}\{|X_i - Y_i|^{\alpha}\mathbb{1}(|X_i - Y_i| > b_n n^{(2-\alpha)/2\alpha})\} \to 0. \qquad \square$$

## Acknowledgements

Many thanks to O. Johnson and R. Samworth for their contributions to the writing of this note. Research partially supported by CNPq, FAPDF, CAPES and FINATEC/UnB.